\documentclass[11pt]{article}
\usepackage{amsfonts,amssymb,amsmath}
-\marginparwidth \oddsidemargin -4mm \evensidemargin \oddsidemargin
\usepackage{amssymb}
\advance\hoffset by 5mm

\newcommand{\eps}{\varepsilon}
\newcommand{\commentout}[1]{}

\newcommand{\la}{\lambda}

\newtheorem{theorem}{Theorem}[section]

\newtheorem{remark}[theorem]{Remark}

\newcommand{\bes}{\begin{displaymath}}
\newcommand{\ees}{\end{displaymath}}
\newcommand{\be}{\begin{equation}}
\newcommand{\ee}{\end{equation}}
\newcommand{\ba}{\begin{eqnarray}}
\newcommand{\ea}{\end{eqnarray}}
\newcommand{\bas}{\begin{eqnarray*}}
\newcommand{\eas}{\end{eqnarray*}}

\newcommand{\B}{{\@Bbb B}}
\newcommand{\C}{{\@Bbb C}}

\newcommand{\F}{{\@Bbb F}}
\renewcommand{\P}{{\@Bbb P}}

\newcommand{\Q}{{\@Bbb Q}}
\newcommand{\bQ}{{\@Bbb Q}}
\newcommand{\N}{{\@Bbb N}}
\newcommand{\R}{{\@Bbb R}}

\newcommand{\W}{{\@Bbb W}}
\newcommand{\Z}{{\mathbb{Z}}}
\newcommand{\al}{\alpha}

\newcommand{\vphi}{\varphi}

\newcommand{\ep}{\varepsilon}

\newcommand{\cA}{\@s A}
\newcommand{\cB}{\@s B}
\newcommand{\cC}{\@s C}
\newcommand{\cD}{\@s D}
\newcommand{\cE}{\@s E}
\newcommand{\cF}{\@s F}
\newcommand{\cG}{\@s G}
\newcommand{\cH}{\@s H}
\newcommand{\cI}{\@s I}
\newcommand{\cJ}{\@s J}
\newcommand{\cK}{\@s K}
\newcommand{\cL}{\@s L}
\newcommand{\cN}{\@s N}
\newcommand{\cM}{\@s M}
\newcommand{\cO}{\@s O}
\newcommand{\cP}{\@s P}
\newcommand{\cR}{\@s R}
\newcommand{\cS}{\@s S}
\newcommand{\cT}{\@s T}
\newcommand{\cV}{\@s V}
\newcommand{\cW}{\@s W}
\newcommand{\cX}{\@s X}
\newcommand{\cY}{\@s Y}
\newcommand{\cZ}{\@s Z}

\newcommand{\bma}{\@bm a}
\newcommand{\bmb}{\@bm b}
\newcommand{\bmc}{\@bm c}
\newcommand{\bmd}{\@bm d}

\newcommand{\bme}{\@bm e}
\newcommand{\bmf}{\@bm f}
\newcommand{\bmg}{\@bm g}
\newcommand{\bmh}{\@bm h}
\newcommand{\bmi}{\@bm i}
\newcommand{\bmj}{\@bm j}
\newcommand{\bmk}{\@bm k}
\newcommand{\bml}{\@bm l}
\newcommand{\bmm}{\@bm m}
\newcommand{\bmn}{\@bm n}
\newcommand{\bmo}{\@bm o}
\newcommand{\bmp}{\@bm p}
\newcommand{\bmq}{\@bm q}
\newcommand{\bmr}{\@bm r}
\newcommand{\bms}{\@bm s}
\newcommand{\bmt}{\@bm t}
\newcommand{\bmu}{\@bm u}
\newcommand{\bmw}{\@bm w}
\newcommand{\bmv}{\@bm v}
\newcommand{\bmx}{\@bm x}
\newcommand{\bx}{\@bm x}
\newcommand{\bmy}{\@bm y}
\newcommand{\bmz}{\@bm z}

\newcommand{\by}{\@bm y}
\newcommand{\bmzero}{\@bm 0}

\newcommand{\gA}{\@g A}
\newcommand{\gD}{\@g D}
\newcommand{\gJ}{\@g J}
\newcommand{\gF}{\@g F}
\newcommand{\gM}{\@g M}
\newcommand{\gR}{\@g R}

\title{On the central limit theorem for some birth and death processes}
 \date{\today}

\author{Tymoteusz Chojecki}

\begin{document}

 \maketitle

 \begin{abstract}
Suppose that $\{X_n,\,n\ge0\}$ is a stationary Markov chain and $V$
is a certain function on a phase space of the chain, called an
observable. We say that the observable satisfies the central limit
theorem (CLT) if   $ Y_n:=N^{-1/2}\sum_{n=0}^NV(X_n)$ converge in
law to a normal random variable, as $N\to+\infty$.  For a stationary
Markov chain with the $L^2$ spectral gap the theorem holds for all
$V$ such that $V(X_0)$ is centered and square integrable, see Gordin
\cite{Go}.  The purpose of this article is  to characterize a family
of observables $V$ for which the CLT holds for a class of birth and
death chains whose dynamics has no spectral gap, so that Gordin's
result cannot be used and the result follows from an application of
Kipnis-Varadhan theory, see \cite{Kip}.
\end{abstract}

\section{Introduction}
 \label{intro}

Suppose that $\{X_n,\,n\ge0\}$ is a stationary Markov chain defined
over a probability space $(\Omega,{\cal F},\mathbb P)$ and $V$ is a
certain function, called an observable, given over the phase of the
chain such that $\mathbb EV(X_0)=0$ and $\mathbb EV^2(X_0)<+\infty$.
Here $\mathbb E$ is the mathematical expectation corresponding to
$\mathbb P$. We say that the chain satisfies the central limit
theorem (CLT) if  the random variables \be
Y_n:=N^{-1/2}\sum_{n=0}^NV(X_n)\label{randomy}\ee converge in law to
a normal random variable, as $N\to+\infty$. Characterization of
Markov chains and the class of observables for which the CLT holds,
is one of the fundamental problems in probability theory. One of the
first results of this type has been the  CLT proved by Doeblin
\cite{Do} for chains whose transition probabilities satisfy what is
now called {\em Doeblin condition}. For stationary Markov chains
with the $L^2$ spectral gap the theorem has been proved by Gordin
\cite{Go}. A remarkable result giving a complete characterization of
reversible Markov chains, i.e. such that for any $N\ge0$ the laws of
$(X_0,X_1\ldots,X_N)$ and of $(X_N,X_{N-1},\ldots,X_0)$ are
identical,  satisfying   the CLT has been proved by Kipnis and
Varadhan in their seminal article \cite{Kip}, see also De Masi  et
al. \cite{Demasi}. It has been shown in \cite{Kip} that the CLT
holds for such chains if only $V$ satisfies
\begin{equation}
\label{eqt-1}
D^2(V):=\lim_{N\to+\infty}\frac{1}{N}\mathbb E\left[\sum_{n=0}^NV(X_n)\right]^2<+\infty.
\end{equation}
One can also prove, see \cite{Kip}, that the limit appearing in
\eqref{eqt-1} always exists,
 being finite or infinite for any reversible chain.
 Sometimes, however, it is not easy to verify the condition directly.

It can be shown that any ergodic and  Markov chain with a finite
phase space has
 a spectral gap, so the CLT  is valid by an application of the Gordin's result.
According to our knowledge, the  examples of chains  not  having the
spectral gap property, yet satisfying the theorem,  concern the
situation when the phase space is
  uncountable,  e.g. tagged particle in a simple exclusion process,  random walks in random environments, diffusions in random media etc., see for instance
  \cite{Demasi,Liggett}. One of the latest review articles about CLT for
  tagged particles and diffusion in random environment is
  \cite{Ola2}.
   The objective of this paper is  to show an application of  the Kipnis-Varadhan theory in the perhaps  simplest possible case (outside  finite chains), namely to a reversible chain with a
   countable phase space but with  no spectral gap (the Gordin's result cannot be used then).  An example like this  is furnished by a
      birth and death chain from Lamperti's problem (see the definition in Section \ref{sec2}), whose phase space  is the
      set of non-negative integers. 
      In the situation considered in the present article
      we give also a necessary and sufficient explicit condition for an observable  $V$ so that \eqref{eqt-1}
      holds. The problem of CLT for the trajectory of the chain has
      been solved in Menshikov and Wade article \cite{Wade}.

As far as  the structure of our paper is concerned,
 in the next section
 we will introduce some basic terms and present three theorems which are our
main results. The  remaining three sections deal with the proofs of
these theorems.

%

\section{Preliminaries and statements of the main results}

\label{sec2}

\subsection{Generalities} Assume that $\{X_n,\ n\geq0\}$ is a Markov chain whose state space
is $\mathbb{Z}_+=\{0,1,2,\ldots\}$. It means that there exists a function $p:\mathbb{Z}_+\times \mathbb{Z}_+\to[0,1]$ such that $\int_{\mathbb{Z}_+} p(x,y)dy=1$ for all  $x\in \mathbb{Z}_+$ and
$$
\mathbb{P}\big[X_{n+1}=x_{n+1}|X_0=x_0,\ldots,X_n=x_n\big]=p(x_n,x_{n+1}).
$$
Here $\int_{\mathbb{Z}_+}f(x)a(dx):=\sum_{0}^{+\infty}f(x)a(x)$ for
an arbitrary $f:\mathbb{Z}_+\to\mathbb R$ and
$a:\mathbb{Z}_+\to[0,+\infty)$. We will write
$\int_{\mathbb{Z}_+}f(x)dx$ when $a(x)\equiv x$. Then
$$
Pf(x)=\int_{\mathbb{Z}_+}p(x,y)f(y)dy,\quad\forall\,f\in B_b(\Z_+)
$$
is called a transition operator.  Here $B_b(\Z_+)$ denotes the space
of all bounded functions on $\Z_+$. Suppose that
$\pi:\mathbb{Z}_+\to(0,1]$ is a strictly positive probability
measure, i.e. $\int_{\mathbb{Z}_+} \pi(x)dx=1$. It is assumed to be
reversible and ergodic with respect to the chain. Ergodicity means
that for any bounded $f$ equality
 $Pf=f$ implies that $f$ is constant  $\pi$ a.s. We say that the
 chain is irreducible if for any $x,y\in\Z_+$ exists $n\geq1$ such
 that $p^n(x,y)>0$, where $p^n(x,y)$ denotes the probability of
 going from $x$ to $y$ in $n$ steps.
\begin{remark}
\em{ It is well known (\cite{Durrett}, p. 338) that a stationary, irreducible Markov chain with a countable state space is ergodic. }
\end{remark}  
Reversibility, on the other hand, means that the  detailed balance
condition holds, i.e.:
 \be
p(x,y)\pi(x)=p(y,x)\pi(y),\quad  \forall\,x,y\in\Z_+. \ee This
condition is equivalent to the fact that $P$ can be extended to a
bounded and symmetric contraction on $L^2(\pi)$ - the space of all
$f$ such that $ ||f||_{\pi}^2
=\int_{\mathbb{Z}_+}f^2(x)\pi(dx)<+\infty$.
 In consequence the spectrum of $P$ lies in $[-1,1].$ Note that
$\lambda_0=1$ is the largest eigenvalue of $P$ corresponding to an
eigenfuction $f_0(x)\equiv 1$. Let
 \be
\lambda_1:=\sup[\langle Pf,f\rangle_{\pi}, \int_{\mathbb{Z}_+}f(x)\pi(dx)=0,\
||f||_{\pi}=1]. \ee We say that the chain has the spectral gap
property, when $\lambda_1<1$. Here $\langle\cdot,\cdot\rangle_{\pi}$ is the scalar product corresponding to $\|\cdot\|_{\pi}$.

Now we will formulate our first main result. It is a simple
criterion for continuity of the spectrum at 1. For any $x\in\Z_+$
let us define $\hat\lambda_0^{(x)}$ by \be
\hat\lambda_0^{(x)}=\sup[\langle Pf,f\rangle_{\pi},\ f(x)=0,
||f||_{\pi}=1].\label{lambda0} \ee We note that
$\hat\lambda_0^{(x)}$ is the largest eigenvalue for the ''reduced''
operator $\hat P=\Pi_xP\Pi_x$, where $\Pi_x$ is the orthogonal
projection onto the subspace $H_x^{(\pi)}:=[f\in L^2(\pi):f(x)=0]$.
\begin{theorem}[A criterion for continuity of the spectrum at
1]$\newline$
\label{thm010810}Suppose that the chain $\{X_n,\,n\ge0\}$ is reversible,
irreducible and there exist $x\in\Z_+$ for which
$\hat\la_0^{(x)}=1$. Then $1$ is not an isolated point of the
spectrum of the transition operator $P$.$\label{TW1}$
\end{theorem}
The proof of  this theorem is presented in Section \ref{sec3}.

\subsection{Birth and death processes} We recall the definition of a birth and death process (see
\cite{Durrett}, p. 295, Example 3.4.). In our setting it is a Markov
chain on countable state space $\Z_+:=\{0,1,2,\ldots,\}$ with the transition
probabilities satisfying that $p(x,y)=0$  iff $|x-y|\not=1$. We also
require that $p_0:=p(0,1)=1$ and
$$
p(x,x+1)=p_x,\qquad p(x,x-1)=q_x,\ x\geq1,
$$
are all strictly positive. Of course we have $p_x+q_x=1$ for all $x\ge1$. In this case there is the
 measure
$ \tilde{\pi}(x)=\prod_{k=1}^{x}(p_{k-1}/q_k) $ which is reversible
and
 unique up to a multiplicative constant (\cite{Durrett} p. 301, Example 4.4). This measure is not necessarily finite. In fact
   \cite{Durrett}, p. 306, Theorem 4.5,
   it is infinite iff the chain is recurrent but not positive recurrent. In case it is positive recurrent
 we have $Z:=\int_{\Z_+}\tilde{\pi}(x)dx<+\infty$ (\cite{Durrett} p. 307, Theorem 4.7) and then $
\pi(x):=Z^{-1}\tilde{\pi}(x) $ is a unique invariant law of the
chain.
 A reader can find more information about this class and more details about
this chain in volume I of Feller's monography \cite{Fe}.

Our second goal is to find relation between transition probabilities
and spectral gap property in the birth and death process.

\begin{theorem}
Let $\{X_n\ ,n\geq0\}$ be a birth and death process with transition
probabilities as above. We have three possible situations then:\\
(i) if $ \lim_{x\to+\infty}p_x=p$, $\lim_{x\to+\infty}q_x=q$ and
$p<q, $ then the chain is positive recurrent and has the spectral
gap
property,\\
(ii) if $p_x=1/2-c_x$ and $q_x=1/2+c_x$,
where \ba
\label{010810} &&
0<c_*=\liminf_{x\to\infty}c_x x^{\al}\le \limsup_{x\to\infty}c_xx^{\al}=c^*<+\infty ,
\ea and
$\al\in(0,1)$, then we have the positive recurrence but we do not
have
the spectral gap property,\\
(iii) if \eqref{010810} holds but  for $\al>1$, then the chain is
recurrent but not positive recurrent.$\label{TW2}$
\end{theorem}
The proof of this result is presented in Section \ref{sec4}.
\begin{remark}
\em{From section I.12, p. 71--76 of \cite{Chung} we know that, when
$\al=1$ then both positive recurrence and null recurrence may occur.
It depends on the constants $c_*,\ c^*$.}
\end{remark}
\begin{remark}
\em{The above result can be interpreted as follows. If we have a
strong drift to the left, i.e. the local drift $D_x:=p_x-q_x$
satisfies $\limsup_{x\to\infty}D_x<0$, then the chain is positive
recurrent and has the spectral gap property. When we have a weaker
drift to the left but $D_x\sim -c/x^{\al}$, for some $c>0$ and
$\al\in(0,1)$, then
 the chain is positive recurrent but does not have the spectral gap.
Finally, when we further increase the probability of
going to the right so that $D_x\sim -c/x^{\al}$, for $\al>1$ and some $c>0$,
 then the chain loses the property of the positive recurrence.}
\end{remark}

The most interesting case of the previous theorem is part (ii). In
this situation we wish to characterize the class of observables for
which the random variables $\eqref{randomy}$ satisfy the CLT The
necessary and sufficient condition for this can be stated as
follows.
\begin{theorem}
Let $V:\Z_+\to\mathbb{R}$ be a zero-mean function in $L^2(\pi)$ and
$\{X_n,\ n\geq0\}$ is the chain from Theorem \ref{TW2} (ii). Then,
$Y_n$ given by $\eqref{randomy}$ satisfies $\eqref{eqt-1}$ iff
$$
\int_{\Z_+}\frac{dx}{\pi(x)}\bigg[\int_{0\leq y\leq
x}V(y)\pi(dy)\bigg]^2<\infty.
$$
 $\label{TW3}$
\end{theorem}

\section{Proof of Theorem \ref{TW1}}

\label{sec3} In order to make our calculations easier, we change our
space $L^2(\pi)$ into space $\ell^2$ corresponding to the counting
measure on $\Z_+$ while $\langle\cdot,\cdot\rangle$, $||\cdot||$
denote the respective scalar product and the $\ell^2$ norm.

We will introduce some terms which are useful throughout the proof. 
Let us denote by $A=[a(x,y)]_{x,y\in\Z_+}$ a matrix with
$a(x,y):=\pi^{1/2}(x)p(x,y)\pi(y)^{-1/2}$. Note that the definition
means that $A=DPD^{-1}$, where
$D=$diag$[\pi^{1/2}(0),\pi^{1/2}(1),\ldots]$, i.e. the operators are
unitary equivalent. In particular the above means that \be
\la_1=\sup[\langle A f, f\rangle;\,\|f\|=1,\,f\in
\ell^2_0].\label{lambda1} \ee Observe that $A$ is a symmetric
matrix. Moreover, for any $f$ with $||f||\le 1$ we have
\begin{eqnarray*}
&&
|\langle Af,f\rangle|:=\left|\int_{\Z_+}\int_{\Z_+}a(x,y)f(x)f(y)dxdy\right|= \left|\int_{\Z_+}
\int_{\Z_+}\pi^{1/2}(x)p(x,y)\pi(y)^{-1/2}f(x)f(y)dxdy\right|\\
&& 
\leq\left\{\int_{\Z_+}\int_{\Z_+}\pi(x)p(x,y)\pi^{-1}(y)f(y)^2dxdy\right\}^{1/2}
\left\{\int_{\Z_+}\int_{\Z_+}p(x,y)f(x)^2dxdy\right\}^{1/2}\\&&\le
||f||^2\le 1.
\end{eqnarray*}
The spectrum of $A$ is also contained in $[-1,1]$, in fact because
of the unitary equivalence, it coincides with the spectrum of
$P$. 
 Note that
$f_*:=(\pi^{1/2}(0),\pi^{1/2}(1),\ldots)$ is an eigenvector that
corresponds to an eigenvalue $\la_0=1$.  
Denote by $\ell^2_0$ the
space consisting of $f\in \ell^2$ such that
$\langle f,f_*\rangle=0.$

Denote by $A'$ matrix $A'=\Pi_xA\Pi_x$. We can easily check that
\be\hat\la_0^{(x)}=\sup[\langle A f, f\rangle;\,\|f\|=1,\,f\in
H_x],\label{lambda0Inna}\ee where $H_x:=[f\in
\ell^2 :f(x)=0]$ and $\hat\la_0^{(x)}$ in
$\eqref{lambda0Inna}$ is defined in $\eqref{lambda0}$.

If the supremum is attained at certain $f^{(0)}$, such that
$\|f^{(0)}\|=1$, $f^{(0)}\in H_x$, then we would have to have
$Af^{(0)}=f^{(0)}$ and that would mean $Pg^{(0)}=g^{(0)}$, where
$g^{(0)}=D^{-1}f^{(0)}$. This, however, would imply $g^{(0)}=c1$ for
some constant $c$, or equivalently
$f^{(0)}=c[\pi^{1/2}(0),\pi^{1/2}(1),\ldots]$. Since $f^{(0)}(x)=0$
we would have $c=0$, which leads to a contradiction. The above means
that $1$ is not in the point spectrum of $A'$ and since it does
belong to the spectrum it must be in its continuous part. We show
that the above implies that \be \la_1=1.\label{4Gwiazdki} \ee
Suppose otherwise, i.e. $\lambda_1<1$. Indeed, suppose that
$f^{(n)}\in H_x$ are such that $\|f^{(n)}\|=1$ and
$$
\langle A f^{(n)}, f^{(n)}\rangle\to 1.
$$
Denote by $Q$ the orthogonal projection onto $\ell^2_0$. We have
$$
f^{(n)}=\alpha_nf_*+Qf^{(n)}
$$
and $\|Qf^{(n)}\|^2=1-\alpha^2_n$. Since $\langle
Af_*,Qf^{(n)}\rangle=0$ we have from $\eqref{4Gwiazdki}$
$$
1\leftarrow\langle A f^{(n)}, f^{(n)}\rangle\leq\alpha_n^2+\lambda_1
(1-\alpha^2_n)\to 1.
$$
This implies $\alpha_n^2\to1$ and in consequence $\|Qf^{(n)}\|\to0$.
Suppose that $\alpha_n\to1$. This yields
$$
\pi^{1/2}(x)\le\|f^{(n)}-f_*\|\to0,
$$
which is impossible. On the other hand, if $\alpha_n\to-1$ we have
$$
\pi^{1/2}(x)\le\|f^{(n)}+f_*\|\to0,
$$
which is again impossible. Hence the conclusion of the theorem
follows.

\section{Proof of Theorem \ref{TW2}}
\label{sec4}
We split our proof into two parts. The first one called, "strong
drift to the left", considers the case (i) from Theorem $\ref{TW2}$,
the second, "weaker drift to the left" deals with the cases (ii) and
(iii) from the theorem. The main point is that in case (i), when the
drift to the left is sufficiently strong the chain has a spectral
gap.
\subsection{Strong drift to the left}
In this case we have $ \lim_{x\to+\infty}p_x=p$,
$\lim_{x\to+\infty}q_x=q$ and $p<q. $ First, we  check whether the
chain is positive recurrent i.e. we verify that
$$
Z:=\int_{\Z_+}\tilde{\pi}(x)dx<\infty,
$$
where
$$
\tilde{\pi}(x)=\prod_{k=1}^{x}\frac{p_{k-1}}{q_k},
$$
see Example 4.4, p. 301 of \cite{Durrett}. From  assumption (i) we
know that for all $\ep>0$, exist $i_0>0$, such that for all $k\geq
k_0$ we have
$$
|p_k-p|<\ep \wedge
|q_k-q|<\ep.
$$
 Now we
can see that
$$
\int_{\Z_+}\tilde{\pi}(x)dx=\int_{0\leq x\leq
k_0}\tilde{\pi}(x)dx+\tilde{\pi}(k_0)\sum_{k=k_0+1}^{\infty}\prod_{j=k_0+1}^k\frac{p_{j-1}}{q_j}\leq
c\sum_{k=k_0+1}^{\infty}\bigg(\frac{p+\ep}{q-\ep}\bigg)^{k-k_0}<\infty
$$
for some $c>0$, provided that $p+\eps<q-\eps$. Hence in this case we
know from \cite{Durrett}, p. 307, Theorem 4.7 that we have a positive
recurrence.

Now we show  the spectral gap property. From
\cite{chen}, see Theorem 1.5, p. 10 case 3, we have
 that $\lambda_1<1$ iff $\delta<\infty$, where
$$
\delta:=\sup_{x\geq1}\int_{0\leq y \leq x-1}dy\left\{[\tilde{\pi}(y)
p_y]^{-1}\int_{x\leq y}\tilde{\pi}(y)dy\right\}.
$$
Observe that
\begin{eqnarray*}
&&
\int_{x\leq y}\tilde{\pi}(y)dy=\tilde{\pi}(x)\int_{x+1\leq k}dk\prod_{l=x+1}^k\frac{p_{l-1}}{q_l}\\
&& \leq \tilde{\pi}(x)\int_{x+1\leq
k}\bigg(\frac{p+\ep}{q-\ep}\bigg)^{k-k_0}dk\le C\tilde{\pi}(x)
\end{eqnarray*}
for some constant $C>0$, provided that $p+\eps<q-\eps$. We can write then
\begin{eqnarray*}
&& \delta\le \frac{C}{p-\eps}\sup_{x\geq1}\int_{0\leq y\leq
x-1}\frac{\tilde{\pi}(x)dy}{\tilde{\pi}(y)}=
\frac{C}{p-\eps}\sup_{x\geq1}\int_{0\leq y\leq x-1}dy\prod_{k=y+1}^{x}\frac{p_{k-1}}{q_k}\\
&& \le
\frac{C}{p-\eps}\times\frac{(p+\ep)(q-\ep)^{-1}}{1-(p+\ep)(q-\ep)^{-1}}<+\infty.
\end{eqnarray*}

\subsection{Weaker drift to the left}
Now we will be interested in the case when $p_x=1/2-c_x, \
q_x=1/2+c_x$,
where \ba &&\liminf_{x\to\infty}c_x x^{\al}=c_*>0,\nonumber\\
&& \limsup_{x\to\infty}c_xx^{\al}=c^*<\infty. \nonumber\ea so we can
find positive constants $K,D_1,D_2$ such that for all $x\geq K$ we
have
$$
\frac{D_1}{x^{\al}}\geq c_x\geq\frac{D_1}{x^{\al}}.
$$
We show that when $\al\in(0,1)$ then we do not have the spectral gap
property but we have  the positive recurrence. Also
we show that when $\al>1$, we do not have the positive recurrence.
Thus, $\al=1$ is a critical exponent where the chain loses the
positive recurrence.

It is easy to see that we can find some
positive constant $c$ for which $p_{k-1}/q_k<1-c/k^{\al}$ for $k>K$.
And for such $c$ we have
$$
-\log\frac{p_{k-1}}{q_k}>-\log(1-\frac{c}{k^{\al}})>\frac{c}{k^{\al}},\qquad
k>K.
$$
Hence, using the integral test for convergence we have
$$
\tilde{\pi}(x)<\exp
\big\{-c\sum_{k=1}^{x}\frac{1}{k^{\al}}\big\}<\tilde{c}\exp(-cx^{1-\al}),
\qquad x>K,
$$
where $\tilde c$ denotes a positive constant. From the
comparison test we see that $\int_{\Z_+}\tilde{\pi}(x)dx<+\infty$ when $\al<1$ and the positive recurrence follows. 

On the other hand we also have some positive constant $c'$, for
which $p_{k-1}/q_k>1-c'/k^{\al}$ for $k>K$. But when $\al>1$, we can
easily check, using again the integral test for convergence, that
$$
\tilde{\pi}(x)>\hat{c}\exp(-c''x^{1-\al}),\qquad x>K,
$$
where $\hat{c},\ c''$ is another positive constant. Then $
\tilde{\pi}(x) $ fails the necessary condition for convergence of
the respective series.

Now we will show that we do not have the spectral gap property when
$\al <1$. To do so we use Theorem \ref{thm010810} and show that $1$
cannot be an isolated point of the spectrum. We choose $x=0$ in
condition \eqref{lambda0}. Denote by $A'$  a symmetric matrix
obtained from $A$ by crossing out the $0$-th column and $0$-th row.
We prove that $\sup_{||f||=1}\langle A'f,f\rangle=1$. Let $f_n$ have
$n^{-1/2}$ on the first
$n$ coordinates and the rest of them vanishes, i.e. $f_n:=[n^{-1/2},n^{-1/2},\ldots,n^{-1/2},0\ldots]$. 
A simple computation shows that
 \bes a(x,x+1) = \sqrt{p_xq_{x+1}}=\left\{
\begin{array}{ll}
\sqrt{\frac{1}{2}+c_1} & \textrm{, $x=0$}\\
\sqrt{\frac{1}{4}+\frac{1}{2}c_{x+1}-\frac{1}{2}c_x-c_xc_{x+1}} &
\textrm{, $x>0$}
\end{array} \right.
\ees

 \bes a(x,x-1) =\sqrt{p_{x-1}q_{x}}= \left\{
\begin{array}{ll}
\sqrt{\frac{1}{2}+c_1} & \textrm{, $x=1$}\\
\sqrt{\frac{1}{4}+\frac{1}{2}c_{x}-\frac{1}{2}c_{x-1}-c_xc_{x-1}} &
\textrm{, $x>1$}
\end{array} \right.
\ees and
$$
\langle
A'f_n,f_n\rangle=2\int_{x\geq1}a(x,x+1)f_n(x)f_n(x+1)dx=\frac{1}{n}\int_{1\leq
x\leq n}\sqrt{1+2c_{x+1}-2c_x-c_xc_{x+1}}dx.
$$
Let $\ep>0$. Then, there exists $ i_0$ such that for  $i>i_0$
$$
|2c_{i+1}-2c_i-c_ic_{i+1}|<\ep.
$$
Hence,
$$
\langle A'f_n,f_n\rangle \geq \sqrt{1-\ep}, \textrm{when}\
n\to\infty.
$$
Since $\ep>0$ was arbitrary we have $\hat\lambda_0^{(0)}=1$ and, by
Theorem \ref{thm010810}, we do not have the spectral gap property.

\section{Proof of Theorem \ref{TW3}}
In this section we assume that the hypothesis of Theorem \ref{TW2}
part (ii) holds. We formulate a sufficient and necessary condition
for an observable $V$ so that \eqref{eqt-1} holds for random variables
given by \eqref{randomy}.

We are going to use the result from \cite{Kip},
see (1.8) p. 3.
According to that result the necessary and sufficient condition for
the validity of \eqref{eqt-1} is that $ ||V||_{-1}^2<\infty, $
where
 \be
||V||_{-1}^2:=\sup_{\varphi\in L^2(\pi)}\bigg\{ 2\langle
V,\varphi\rangle_{\pi}-\langle
(I-P)\varphi,\varphi\rangle_{\pi}\bigg\}.\label{def: norma}
 \ee
%
In what follows we
will find the maximizer of this supremum by solving the
Euler-Lagrange equation $\eqref{euler-lagrange}$. The maximizer belongs to a certain Hilbert space, that we denote by  ${\cal H}_1$, which is
a bigger space than $L^2(\pi)$.

Observe that we have a following equality
 \be\langle
(I-P)\varphi,\varphi\rangle_{\pi}=\frac{1}{2}\int_{\Z_+}\int_{\Z_+}p(x,y)(\varphi(y)-\varphi(x))^2dy\pi(dx).
\label{def:rownosc} \ee Denote the discrete gradient $\partial
f(x)=f(x+1)-f(x)$ and its dual, with respect to the scalar product
from  $\ell^2$, $\partial^* f(x)=f(x-1)-f(x)$. It is easy to see
that when the supp$f$ is compact then the following integration by
parts formula holds
$$
\int_{\Z_+}\partial
f(x)g(x)dx=\int_{x\geq1}f(x)\partial^*g(x)dx-f(0)g(0).
$$
In our case we have \be {\cal E}(\vphi):=\langle
(I-P)\varphi,\varphi\rangle_{\pi}=\frac{1}{2}\int_{x\geq1}\bigg[\bigg(\frac{1}{2}+c_x\bigg)\bigg(\partial^*\vphi(x)\bigg)^2
+
\bigg(\frac{1}{2}-c_x\bigg)\bigg(\partial\vphi(x)\bigg)^2\bigg]\pi(dx).\nonumber
\ee We can  find some positive constants $c_1,c_2$, for which
 \ba
 c_1{\cal E}_0(\vphi)\le  {\cal E}(\vphi)\leq
c_2{\cal E}_0(\vphi)\label{def:
e0}, \ea
where
$$
{\cal E}_0(\vphi):=\int_{\Z_+}\big(\partial\vphi(x)\big)^2\pi(dx).
$$
Then, \begin{eqnarray}
\label{020810} &&\sup_{\vphi\in L^2(\pi)}\big\{2\langle V,\vphi\rangle_{\pi}- c_2{\cal E}_0(\vphi)\big\}\le ||V||_{-1}^2\leq\sup_{\vphi\in
L^2(\pi)}\big\{2\langle V,\vphi\rangle_{\pi} -
c_1{\cal E}_0(\vphi)\big\}.
\end{eqnarray}
We can see therefore  that $\|V\|_{-1}<+\infty$ iff the supremum of the functional $\Phi(\varphi)$ appearing on the right hand side of
\eqref{020810} is finite. 
Since the functional $\Phi(\cdot)$ is weakly upper  semicontinuous on a Hilbert space
$$
{\cal
H}_1:=[\varphi:{\cal E}_0(\varphi)<+\infty],
$$
 and $\lim_{\|\varphi\|_{{\cal H}_1}\to+\infty}\Phi(\varphi)=-\infty$, it attains its maximum
\begin{equation}
\label{min} \Phi_*=\sup[\Phi(\varphi):\varphi\in L^2(\pi)]<+\infty
\end{equation}
and
its maximizer $\varphi_*$  has to satisfy the Euler-Lagrange equation, which in this case reads
\begin{equation}
\label{euler-lagrange}
V(x)\pi(x)=\partial^*[\pi(x)\partial\varphi_*(x)]
,\quad\forall\,x\ge1
\end{equation}
and $V(0)=-\partial\varphi(0)$,
or equivalently$$
\pi(x)\partial\varphi(x)=-\pi(x)
V(x)+\pi(x-1)\partial\varphi(x-1),\quad\forall\,x\ge 1
$$
and $\partial\varphi(0)=-V(0)$. Note that from this equation we get
\begin{eqnarray*}
&&\pi(x)\partial\varphi(x) =-\int_{1\leq y\leq
x}V(y)\pi(dy)+\pi(0)\partial\varphi(0)=-\int_{0\leq y\leq
x}V(y)\pi(dy),
\end{eqnarray*}
hence
\begin{eqnarray*}
&&\partial\varphi(x) =-\frac{1}{\pi(x)}\int_{0\leq y\leq
x}V(y)\pi(dy),
\end{eqnarray*}
The supremum in \eqref{min} equals
$$
\Phi_*=\int_{\Z_+}[\partial\varphi(x)]^2\pi(dx)=\int_{\Z_+}\frac{dx}{\pi(x)}\left[\int_{0\leq
y\leq x}V(y)\pi(dy)\right]^2
$$
and the requirement that $\Phi_*<+\infty$ is  the same requirement
as  $\|V\|_{-1}<+\infty$, thus  the conclusion of Theorem \ref{TW3} follows. Note that the fact that $\Phi_*<+\infty$  in particular implies that
$$
\int_{\Z_+}V(y)\pi(dy)=0.
$$
\subsection{Acknowledgement}
I would like to express my deep gratitude to Professor Tomasz
Komorowski whose help and comments have been invaluable and without
whom writing this paper would have been impossible, and to Professor
Andrew Wade, who pointed out my error in part (iii) of Theorem
\ref{TW2}.


\begin{thebibliography}{99}
\bibitem{chen}Chen, M., \emph{Eigenvalues, Inequalities, and Ergodic Theory}, Springer,
(2005).\\
\bibitem{Chung}Chung, K.L., \emph{Markov Chains with Stationary Transition Probabilities}, 2nd edition, Springer-Verlag, Berlin, (1967).

\bibitem{Demasi}De Masi, A. and Ferrari, P. A. and Goldstein, S. and Wick, W. D., \emph{An invariance principle for reversible Markov processes. Applications to random motions in random environments}, J. Statist. Phys., 55,
(1989).\\
\bibitem{Do}Doeblin, W., \emph{Sur deux proble`mes de M. Kolmogoroff concernant les chaines dénombrables}, Bull. Soc. Math. France vol. 66 (1938) pp.
210-220.\\
\bibitem{Durrett} Durrett, R., \emph{Probability theory and examples},
(Wadsworth Publishing Company, Belmont, 1996).\\
\bibitem{Fe}Feller, W. \emph{An introduction to probability theory and its applications}, Vol. II. Second edition John Wiley, Sons, Inc., New York-London-Sydney
(1971).\\
\bibitem{Go}Gordin, M. I., \emph{The central limit theorem for stationary
 processes}, Dokl. Akad. Nauk SSSR,188, pp 739Ð741, (1969).

\bibitem{Kip} Kipnis, C., Varadhan S. R. S. \emph{Central limit
theorem for additive functionals of reversible Markov processes and
applications to simple exclusions}, Comm. Math. Phys. 104 (1986),
no. 1, 1--19.

\bibitem{Liggett} Liggett, T., {\em Stochastic Interacting Systems: Contact, Voter and Exclusion Processes,} Grund. fur Math. Wissen., Vol 324, Springer
(1999).
\bibitem{Wade} Menshikov M., Wade A.,\emph{Rate of escape and central limit theorem for the
supercritical Lamperti problem}, Stochastic Processes and their
Applications 120 (2010) 2078–2099.
\bibitem{Ola2} Olla, S., \emph{Notes on Central Limits Theorems for Tagged Particles and Diffusions in Random
Enviornment}, Etats de la recherche: Milieux Aleatoires CIRM, Luminy
(2000).
\end{thebibliography}
\end{document}